\documentclass{article}
\usepackage{eurosym}
\usepackage{amssymb}
\usepackage{amsfonts}
\usepackage{amsmath}

\begin{document}

\begin{center}
{\Large Some remarks regarding \textit{l-}elements defined in algebras
obtained by the Cayley-Dickson process}

\begin{equation*}
\end{equation*}%
Cristina FLAUT and Diana SAVIN 
\begin{equation*}
\end{equation*}
\end{center}

\textbf{Abstract. }{\small In this paper, we define a special class of
elements in the algebras obtained by the Cayley-Dickson process, called }$l-$%
{\small elements. We find conditions such that these elements to be
invertible.~These conditions can be very useful for finding new identities,
identities which can help us in the study of the properties of these
algebras.} 
\begin{equation*}
\end{equation*}

\textbf{Key Words}: quaternion algebras; octonion algebras; Cayley-Dickson
algebras; special sequences.

\medskip

\textbf{2000 AMS Subject Classification}: 15A24, 15A06, 16G30, 11R52, 11B39,
11R54.%
\begin{equation*}
\end{equation*}

\bigskip

\textbf{1. Introduction}%
\begin{equation*}
\end{equation*}

In the last years, many papers were devoted to the study of some special
sequences or special elements, especially in their connections with
particular cases of algebras obtained by the Cayley-Dickson (quaternions,
octonions). Some of them are: [Ba, Pr; 09], [Ca; 16], [Ci, Ip; 16], [Fa, Pl;
07(1)], [Fa, Pl; 07(2)], [Fl; 09], [Fl, Sa; 15(1)], [Fl, St; 09], [Ha, 15],
[Mo; 17], [Si,18].

Since the algebras obtained by the Cayley-Dickson process, denoted $A_{t}$,
are poor in properties when their dimension increase, losing commutativity,
associativity and alternativity, in this paper, we will find elements with
supplementary properties, namely a special set of invertible elements in $%
A_{t},$ which can help us in the study of these algebras. Algebras obtained
by the Cayley-Dickson process will be presented in the following.

Let $K$ be a commutative field with $char\left( K\right) \neq 2.~$Let $A$ \
be a finite dimensional unitary algebra over a field $\ K.$ We define a 
\textit{scalar} \textit{involution} $\,$%
\begin{equation*}
\,\,\,\overline{\phantom{x}}:A\rightarrow A,a\rightarrow \overline{a},
\end{equation*}%
to be a linear map with the following properties$\,\,\,$%
\begin{equation*}
\overline{ab}=\overline{b}\overline{a},\,\overline{\overline{a}}=a,
\end{equation*}%
$\,\,$and 
\begin{equation*}
a+\overline{a},a\overline{a}\in K\cdot 1\ \text{for all }a,b\in A.\text{ }
\end{equation*}%
An element $\,\overline{a}\in A$ is called the \textit{conjugate} of the
element $a\in A,$ the linear form$\,\,$%
\begin{equation*}
\,\,\mathbf{t}:A\rightarrow K\,,\,\,\mathbf{t}\left( a\right) =a+\overline{a}
\end{equation*}%
and the quadratic form 
\begin{equation*}
\mathbf{n}:A\rightarrow K,\,\,\mathbf{n}\left( a\right) =a\overline{a}\ 
\end{equation*}%
are called the \textit{trace} and the \textit{norm \ }of \ the element $a.$

We consider$\,\,\,\gamma \in K$ a fixed non-zero element. We define the
following algebra multiplication on the vector space $A\oplus A:$ 
\begin{equation}
\left( a_{1},a_{2}\right) \left( b_{1},b_{2}\right) =\left(
a_{1}b_{1}+\gamma \overline{b}_{2}a_{2},a_{2}\overline{b_{1}}%
+b_{2}a_{1}\right) .  \tag{1.1.}
\end{equation}%
\newline
The obtained algebra structure over $A\oplus A,$ denoted by $\left( A,\gamma
\right) ,$ is called the \textit{algebra obtained from }$A$\textit{\ by the
Cayley-Dickson process.} $\,$We have $\dim \left( A,\gamma \right) =2\dim A$.

Let $x\in \left( A,\gamma \right) $, $x=\left( a_{1},a_{2}\right) $. The map 
\begin{equation*}
\,\,\,\overline{\phantom{x}}:\left( A,\gamma \right) \rightarrow \left(
A,\gamma \right) \,,\,\,x\rightarrow \bar{x}\,=\left( \overline{a}_{1},\text{%
-}a_{2}\right) ,
\end{equation*}%
\newline
is a scalar involution of the algebra $\left( A,\gamma \right) $, extending
the involution $\overline{\phantom{x}}\,\,\,$of the algebra $A.$ We take 
\begin{equation*}
\,\mathbf{t}\left( x\right) =\mathbf{t}(a_{1})
\end{equation*}%
and$\,\,\,$ 
\begin{equation*}
\mathbf{n}\left( x\right) =\mathbf{n}\left( a_{1}\right) -\gamma \mathbf{n}%
(a_{2})
\end{equation*}%
to be $\,$the \textit{trace} and the \textit{norm} of the element $x\in $ $%
\left( A,\gamma \right) $.\thinspace $\,$

\thinspace If we consider $A=K$ \thinspace and we apply this process $t$
times, $t\geq 1,\,\,$we obtain an algebra over $K,\,\,$denoted by%
\begin{equation}
A_{t}=\left( \frac{\gamma _{1},...,\gamma _{t}}{K}\right) .  \tag{1.2. }
\end{equation}

Using induction, in this algebra, the set $\{1,e_{2},...,e_{n}\},n=2^{t},$
generates a basis with the properties:%
\begin{equation}
e_{i}^{2}=\gamma _{i}1,\,\,_{i}\in K,\gamma _{i}\neq 0,\,\,i=2,...,n 
\tag{1.3.}
\end{equation}%
and \ 
\begin{equation}
e_{i}e_{j}=-e_{j}e_{i}=\beta _{ij}e_{k},\,\,\beta _{ij}\in K,\,\,\beta
_{ij}\neq 0,i\neq j,i,j=\,\,2,...n,  \tag{1.4.}
\end{equation}%
$\ \ \beta _{ij}$ and $e_{k}$ being uniquely determined by $e_{i}$ and $%
e_{j}.$ We remark that if $x\in A_{t}$ and $n\left( x\right) \neq 0,$ then $%
x $ is an invertible element in $A_{t}.$

A finite-dimensional algebra $A$ is \textit{a division} algebra if and only
if $A$ does not contain zero divisors. For other details regarding algebras
obtained by the Cayley-Dickson process, the reader is referred to [Sc;66].

For $t=2,$ we obtain the generalized quaternion algebras, denoted $\,\mathbb{%
H}(\gamma _{1},\gamma _{2}),$ and for $t=3,~$we obtain the generalized
octonion algebras, denoted $\mathbb{O}(\gamma _{1},\gamma _{2},\gamma _{3})$.

For example, for $x\in \mathbb{H}(\gamma _{1},\gamma _{2}),$ $%
x=x_{0}+x_{1}e_{1}+x_{2}e_{2}+x_{3}e_{3},$ its norm is 
\begin{equation}
\,\mathbf{n}\left( x\right) =x\overline{x}=x_{0}^{2}-\gamma
_{1}x_{1}^{2}-\gamma _{2}x_{2}^{2}+\gamma _{1}\gamma _{2}x_{3}^{2}\in K. 
\tag{1.5.}
\end{equation}

If, \thinspace for $x\in \mathbb{A},$ where $\mathbb{A=\{\,\mathbb{H}(}%
\gamma _{1},\gamma _{2}\mathbb{),\mathbb{O}(}\gamma _{1},\gamma _{2},\gamma
_{3}\mathbb{)\}},$ the relation $\mathbf{n}\left( x\right) =0$ implies $x=0$%
, then these algebras are division algebras. A quaternion and an octonion
non-division algebras are called \textit{split} algebras. Using the above
notations, we remark that the quaternion algebra $\mathbb{H}\left(
-1,-1\right) =\left( \frac{-1,-1}{\mathbb{R}}\right) $ and the octonion
algebra $\mathbb{O}\left( -1,-1,-1\right) =\left( \frac{-1,-1,-1}{\mathbb{R}}%
\right) $ are division algebras. If the field $K$ is a finite field, then
the quaternion and octonion algebras are always split algebras.

An algebra $A\in \mathbb{A}$ such that $n\left( xy\right) =n\left( x\right)
n\left( y\right) ,$ for all $x,y\in A,$ is called a \textit{composition
algebra}. For $t\geq 4,~$the algebras $A_{t}$ are not composition algebras.

A unitary algebra $A\neq K$ such that we have $x^{2}+\alpha _{x}x+\beta
_{x}=0,$ for each $x\in A,$ with $\alpha _{x},\beta _{x}\in K,$ is called a 
\textit{quadratic algebra}. An algebra $A$ with a scalar involution is
quadratic. From here, we get that the algebras obtained by the
Cayley-Dickson process are quadratic.

An algebra $A$ is called \textit{alternative} if $x^{2}y=x\left( xy\right) $
and $xy^{2}=\left( xy\right) y,$ for all $x,y\in A,$ \textit{\ flexible} if $%
x\left( yx\right) =\left( xy\right) x=xyx,$ for all $x,y\in A$ and \textit{%
power associative} if the subalgebra $<x>$ of $A,$ generated by any element $%
x\in A,$ is associative.

Algebras $A_{t}$ are noncommutative, for $t\geq 2,$ nonassocitive, for $%
t\geq 3,$ non-alternative, for $t\geq 4$ and, in general, non division
algebras. They are power associative and flexible for all $t.$

\bigskip

\begin{equation*}
\end{equation*}

\textbf{2.} $l$\textbf{-numbers and some of their properties} 
\begin{equation*}
\end{equation*}

In [Ho; 61] and [Fl, Sa; 17], were defined and studied the $%
(a,b,x_{0},x_{1}) $ numbers. If we consider $l$ a nonzero natural number,
for $a=l,b=1,x_{0}=0,x_{1}=1,~$in [Sa; 18] was considered the following
sequence, where 
\begin{equation}
p_{n}=l\cdot p_{n-1}+p_{n-2},\;n\geq 2,p_{0}=0,p_{1}=1.  \tag{2.1.}
\end{equation}%
Since these numbers are $\left( l,1,0,1\right) $ numbers, we will call them $%
l-$\textit{numbers}. For $l=1,$ it is obtained the Fibonacci numbers and for 
$l=2,$ it is obtained the Pell numbers.\medskip \newline

\textbf{Remark 2.1.} ([Sa; 18]). \ Let $\left( p_{n}\right) _{n\geq 0}$ be
the sequence previously defined. Then, the following relations are true:%
\newline
1)%
\begin{equation*}
p_{n}^{2}+p_{n+1}^{2}=p_{2n+1},\ \left( \forall \right) n\in \mathbb{N}.
\end{equation*}%
2) For $\alpha =\frac{l+\sqrt{l^{2}+4}}{2}$ and $\beta =\frac{l-\sqrt{l^{2}+4%
}}{2},$ we obtain that \smallskip\ 
\begin{equation*}
p_{n}=\frac{\alpha ^{n}-\beta ^{n}}{\alpha -\beta }=\frac{\alpha ^{n}-\beta
^{n}}{\sqrt{l^{2}+4}},\ \ \left( \forall \right) n\in \mathbb{N},
\end{equation*}%
called the \textbf{Binet's formula for the sequence} $\left( p_{n}\right)
_{n\geq 0}.\medskip $

\textbf{Proposition 2.2.} \textit{With the above notations, the following
relations hold:}

1) \textit{If} $d\mid n,$ \textit{then } $p_{d}\mid p_{n}.$

2) $p_{m+n}=p_{m}p_{n+1}+p_{m-1}p_{n}.\medskip $

\bigskip \textbf{Proof.} 1) We use the Binet formula from the above remark.
If $d\mid n,n=dc,$ 
\begin{equation*}
p_{n}=\frac{\alpha ^{dc}-\beta ^{dc}}{\sqrt{l^{2}+4}}=\frac{(\alpha
^{d})^{c}-(\beta ^{d})^{c}}{\sqrt{l^{2}+4}}=
\end{equation*}%
\begin{equation*}
=\frac{(\alpha ^{d}-\beta ^{d})M}{\sqrt{l^{2}+4}},
\end{equation*}%
therefore $p_{d}\mid p_{n}.$

2) We use induction after $n.$ For $n=0,$ we have $p_{m}=p_{m}.$ Supposing
affirmation true for each $k$$\leq$ $n,$ we will prove for $n+1.$ We have $%
p_{m+n+1}=lp_{m+n}+p_{m+n-1}=l\left( p_{m}p_{n+1}+p_{m-1}p_{n}\right) +$%
\newline
$+p_{m}p_{n}+p_{m-1}p_{n-1}=p_{m}\left( lp_{n+1}+p_{n}\right) +p_{m-1}\left(
lp_{n}+p_{n-1}\right) =$\newline
$=p_{m}p_{n+2}+p_{m-1}p_{n+1}.\Box \smallskip $

\medskip

\textbf{Proposition 2.3.} \textit{Let} $(p_{n})_{n\geq 0}$ \textit{be the} $%
l-$\textit{sequence defined in }$\left( 2.1\right) $. \textit{Then, the
following relations are true:}\newline
i) 
\begin{equation}
p_{n}+p_{n+4}=\left( l^{2}+2\right) \cdot p_{n+2};  \tag{2.2.}
\end{equation}%
ii) 
\begin{equation}
p_{n}+p_{n+8}=\left[ \left( l^{2}+2\right) ^{2}-2\right] \cdot p_{n+4}; 
\tag{2.3.}
\end{equation}%
\textit{For} $k\geq 3,$ \textit{we have} 
\begin{equation}
p_{n}+p_{n+2^{k}}=\underbrace{\left[ \left( \left( \left( l^{2}+2\right)
^{2}-2\right) ^{2}...-2\right) ^{2}-2\right] }_{k-2\;times\;of\;-2}\cdot
p_{n+2^{k-1}}.  \tag{2.4.}
\end{equation}

\textbf{Proof.} i) With the above notations, we remark that $\alpha ^{2}+%
\frac{1}{\alpha ^{2}}=\alpha ^{2}+\beta ^{2}=\left( \alpha +\beta \right)
^{2}-2\alpha \beta =l^{2}+2.~$\ Applying Binet's formula for $l-$ sequence,
we have: 
\begin{equation*}
p_{n}+p_{n+4}=\frac{\alpha ^{n}-\beta ^{n}}{\alpha -\beta }+\frac{\alpha
^{n+4}-\beta ^{n+4}}{\alpha -\beta }=
\end{equation*}%
\begin{equation*}
\frac{\alpha ^{n+2}\cdot \left( \alpha ^{2}+\frac{1}{\alpha ^{2}}\right)
-\beta ^{n+2}\cdot \left( \beta ^{2}+\frac{1}{\beta ^{2}}\right) }{\alpha
-\beta }=
\end{equation*}%
\begin{equation*}
\left( l^{2}+2\right) \cdot \frac{\alpha ^{n+2}-\beta ^{n+2}}{\alpha -\beta }%
=\left( l^{2}+2\right) p_{n+2}.
\end{equation*}%
ii) Using i), we obtain: 
\begin{equation*}
p_{n}+p_{n+8}=\left( p_{n}+p_{n+4}\right) +\left( p_{n+4}+p_{n+8}\right)
-2\cdot p_{n+4}=
\end{equation*}%
\begin{equation*}
=\left( l^{2}+2\right) \cdot p_{n+2}+\left( l^{2}+2\right) \cdot
p_{n+6}-2\cdot p_{n+4}=
\end{equation*}%
\begin{equation*}
=\left( l^{2}+2\right) \cdot \left( p_{n+2}+p_{n+6}\right) -2\cdot p_{n+4}=
\end{equation*}%
\begin{equation*}
=\left[ \left( l^{2}+2\right) ^{2}-2\right] \cdot p_{n+4}.
\end{equation*}%
iii) We use induction after $k$$\in $$\mathbb{N},$ $k$$\geq 3$$,$ to prove
the following statement: 
\begin{equation*}
P\left( k\right) :\ p_{n}+p_{n+2^{k}}=\underbrace{\left[ \left( \left(
\left( l^{2}+2\right) ^{2}-2\right) ^{2}...-2\right) ^{2}-2\right] }%
_{k-2\;times\;of\;-2}\cdot p_{n+2^{k-1}}.
\end{equation*}%
In i) we proved $P\left( 2\right) ,$ in ii) we proved $P\left( 3\right) .$
Supposing that $P\left( k\right) $ is true, it is easy to prove $P\left(
k+1\right) ,$ in the same way as we passed from $P\left( 2\right) $ to $%
P\left( 3\right) .\Box \smallskip $

\textbf{Remark 2.4.} We will denote%
\begin{equation}
M_{k}=\underbrace{\left[ \left( \left( \left( l^{2}+2\right) ^{2}-2\right)
^{2}...-2\right) ^{2}-2\right] }_{k-2\;times\;of\;-2},k\geq 3.  \tag{2.5.}
\end{equation}%
We remark that $M_{k+1}=M_{k}^{2}-2.$ We put $M_{2}=$ $l^{2}+2.$ It results, 
$M_{k}>0,$ for all $k\in \mathbb{N}$, $k\geq 2$.

\begin{equation*}
\end{equation*}

\textbf{3. }$l$-\textbf{elements in algebras obtained by the Cayley-Dickson
process}

\begin{equation*}
\end{equation*}

With the above notations, we will denote 
\begin{equation*}
S_{t}=M_{2}M_{3}M_{4}...M_{t-2}M_{t-1}M_{t},t\geq 2.
\end{equation*}%
We will put $S_{1}=1.$

Let $A_{t}=\left( \frac{\gamma _{1},...,\gamma _{t}}{K}\right) $ be an
algebra obtained by the Cayley-Dickson process and $n$$\in $$\mathbb{N}.$ We
define the $n$-th $l-$element $P_{n}\in A_{t}~$to be an element of the form 
\begin{equation}
P_{n}\text{=}p_{n}\cdot 1+p_{n+1}\cdot e_{1}+p_{n+2}\cdot e_{2}+p_{n+3}\cdot
e_{3}+...+p_{n+2^{t}-1}e_{t-1}.  \tag{3.1.}
\end{equation}

\textbf{Proposition 3.1.} \textit{Let} $v$ \textit{be a real number}, $q$ 
\textit{be a prime number, }$K\in \{\mathbb{Q},\mathbb{R,Z}_{q}\},~$\textit{%
and} $A_{t}=\left( \frac{-1,...,-1,v}{K}\right) $ \textit{be an algebra
obtained by the} \textit{Cayley-Dickson process of dimension} $r=2^{t}$%
\smallskip , $t\geq 2,$\ \textit{and} 
\begin{equation*}
P_{n}=p_{n}\cdot 1+p_{n+1}\cdot e_{1}+p_{n+2}\cdot e_{2}+p_{n+3}\cdot
e_{3}+...+p_{n+2^{t}-1}e_{t-1},
\end{equation*}%
\textit{be an }$l-$\textit{element. Therefore, the following statements are
true:}

1) \textit{We have} 
\begin{equation}
\boldsymbol{n}\left( P_{n}\right)
=S_{t-1}(p_{2n+2^{t-1}-1}-vp_{2n+2^{t+1}-1-2^{t-1}}).  \tag{3.2.}
\end{equation}

2) \textit{If }$v\in \mathbb{R}-\left( -1,1\right) ,$ t\textit{he} $n$-%
\textit{th }$l-$\textit{element \ }$P_{n}$ \textit{is invertible in }$%
A_{t}=\left( \frac{-1,...,-1,v}{K}\right) ,$ \textit{for all} $n$$\in $$%
\mathbb{N}$.

3) \textit{If }$v=-1,$ \textit{we have} 
\begin{equation}
\boldsymbol{n}\left( P_{n}\right) =S_{t}p_{2n+2^{t}-1}.  \tag{3.3.}
\end{equation}

\textbf{Proof.} 1) Using Remark 2.1 and Proposition 2.3, we compute the
norm. In this case, we obtain: 
\begin{equation*}
\boldsymbol{n}\left( P_{n}\right) \text{=}p_{n}^{2}\text{+}p_{n+1}^{2}\text{+%
}p_{n+2}^{2}\text{+}p_{n+3}^{2}\text{+...+}p_{n+2^{t-1}-1}^{2}\text{-}%
v(p_{n+2^{t-1}}^{2}\text{+...+}p_{n+2^{t}-2}^{2}\text{+}p_{n+2^{t}-1}^{2})%
\text{=}
\end{equation*}%
\begin{equation*}
=p_{2n+1}+p_{2n+5}+p_{2n+9}+p_{2n+13}+...-v(p_{2n+2^{t}+1}+...+p_{2n+2^{t+1}-3})=
\end{equation*}%
\begin{equation*}
=\left( l^{2}+2\right) p_{2n+3}+\left( l^{2}+2\right)
p_{2n+11}+...-v(...+\left( l^{2}+2\right) p_{2n+2^{t+1}-5})=
\end{equation*}%
\begin{equation*}
\text{=}\left( l^{2}\text{+}2\right) \text{\{[}\left( l^{2}\text{+}2\right)
^{2}\text{-}2\text{]}p_{2n\text{+}7}\text{+[}\left( l^{2}\text{+}2\right)
^{2}\text{-}2\text{]}p_{2n\text{+}23}\text{+...-}v\text{(...+[}\left( l^{2}%
\text{+}2\right) ^{2}\text{-}2\text{]}p_{2n\text{+}2^{t+1}-9}\text{)\}=}
\end{equation*}%
\begin{equation*}
M_{2}M_{3}[M_{4}p_{2n+15}+...-v(...+M_{4}p_{2n+2^{t+1}-17})]=
\end{equation*}%
\begin{equation*}
.....
\end{equation*}

\begin{equation*}
=M_{2}M_{3}M_{4}...M_{t-2}[M_{t-1}p_{2n+2^{t-1}-1}-vM_{t-1}p_{2n+2^{t+1}-1-2^{t-1}})]=
\end{equation*}

\begin{equation*}
=M_{2}M_{3}M_{4}...M_{t-2}M_{t-1}(p_{2n+2^{t-1}-1}-vp_{2n+2^{t+1}-1-2^{t-1}}).
\end{equation*}

2) Indeed, since $2n+2^{t-1}-1<2n+2^{t+1}-1-2^{t-1},$ we have that $%
p_{2n+2^{t+1}-1-2^{t-1}}>p_{2n+2^{t-1}-1},$ since $(p_{n})_{n\geq 0},$ is an
increasing sequence. From here, we get that $\boldsymbol{n}\left(
P_{n}\right) \neq 0,$ even if $v$ is positive or negative$.$ Therefore $%
P_{n} $ is an invertible element in $A_{t}=\left( \frac{-1,...,-1,v}{K}%
\right) .$

3) Taking $v=-1$ in relation (3.2), we get 
\begin{equation*}
\boldsymbol{n}\left( P_{n}\right)
=M_{2}M_{3}M_{4}...M_{t-2}M_{t-1}M_{t}(p_{2n+2^{t-1}-1}+p_{2n+2^{t+1}-1-2^{t-1}})=
\end{equation*}%
\begin{equation*}
=M_{2}M_{3}M_{4}...M_{t-2}M_{t-1}M_{t}p_{2n+2^{t}-1}.
\end{equation*}%
$\Box \smallskip $\smallskip \qquad

An idempotent element in a unitary ring, with $0\neq 1,$ is an element $x$, $%
x\neq 0$ and $x\neq 1,$ such that $x^{2}=x.\medskip \bigskip $

\textbf{Theorem 3.2. } \textit{With the above notations, let} \smallskip\ $q$
\textit{be a prime number and let} $d$ \ \textit{be the first positive
integer such that} $q\mid p_{d}.$

1)\textit{\ Therefore} $q\mid p_{n}$ \textit{if and only if} $d\mid n.$

2) \textit{If} $\gcd \left( q,S_{t}\right) =1,~t\geq 2,$ \textit{and} $d$ 
\textit{is an even number, then all} $l-$\textit{elements }$P_{n},n\geq d,$ 
\textit{from the algebra} $A_{t}=\left( \frac{-1,...,-1}{\mathbb{Z}_{q}}%
\right) $ \textit{are invertible}.

3) \textit{If} $t\in \{2,3\},$\ $\gcd \left( q,S_{t}\right) =1,~t\geq 2,$ 
\textit{and}\ $d$ \textit{is an odd number such that} $d\mid (2^{t}-1),$ 
\textit{then} \textit{there are no idempotent }$l-$\textit{elements in the
algebra }$A_{t}=\left( \frac{-1,...,-1}{\mathbb{Z}_{q}}\right) .$

\textbf{Proof. }

1) Supposing that $q\mid p_{n},$ let $n=dc+r,0\leq r<d,$ therefore 
\begin{equation*}
p_{n}=p_{dc+r}=p_{r}p_{dc+1}+p_{r-1}p_{dc},
\end{equation*}%
from Proposition 2.2$.$ Since $q\mid p_{d},$ we have that $q\mid p_{dc}.$
Using that $q\mid p_{n},$ it results $q\mid p_{r}p_{dc+1}.$ If $q\mid p_{r}$%
, we have a contradiction with the choice of the number $d.$ If $q\mid
p_{dc+1},$ we have \ $q\mid p_{dc-1}$ and so on. Then we obtain that $q\mid
p_{0},$ false. Therefore, $r=0$ and $d\mid n.$ Conversely, assuming that $%
d\mid n,$ we have $n=dc,~c$ a natural number. From here, we have $p_{d}\mid
p_{n},$ then $q\mid p_{n}.$

2) Indeed, since $\boldsymbol{n}\left( P_{n}\right)
=S_{t}p_{2n+2^{t}-1}\equiv 0($mod$~q),$ we obtain $p_{2n+2^{t}-1}\equiv 0$
(mod $q$). It results $d/\left( 2n+2^{t}-1\right) ,$ which is false.

3) Supposing that there are $l-$elements $P_{n}$$\in A_{t}=\left( \frac{%
-1,...,-1}{\mathbb{Z}_{q}}\right) $ such that $P_{n}^{2}=P_{n}$ in $A_{t},$
we obtain $P_{n}\cdot \left( P_{n}-1\right) =0.$ Since $t\in \{2,3\},$ the
algebra $A_{t}$ is a composition algebra. It results that $P_{n}$ and $%
P_{n}-1$ are zero divisors in the algebra $A_{t}$ and their norm are zero.
These are equivalent with the system

\begin{equation*}
\left\{ 
\begin{array}{c}
\boldsymbol{n}\left( P_{n}\right) =\widehat{0} \\ 
\boldsymbol{n}\left( P_{n}-1\right) =\widehat{0}%
\end{array}%
\right.
\end{equation*}%
in $\mathbb{Z}_{q}.$ Using Proposition 3.1, we get

\begin{equation*}
\left\{ 
\begin{array}{c}
S_{t}p_{2n+2^{t}-1}\equiv 0(\text{mod}~q) \\ 
(p_{n}-1)^{2}+p_{n+1}^{2}+p_{n+2}^{2}+...+p_{n+2^{t}-1}^{2}\equiv 0(\text{mod%
}~q).%
\end{array}%
\right. 
\end{equation*}%
Since $\gcd \left( q,S_{t}\right) =1,$ we obtain $p_{2n+2^{t}-1}=0($mod$~q),$
that means $q\mid p_{2n+2^{t}-1}.$ It results, from above, that 
\begin{equation*}
2n+2^{t}-1\equiv 0(\text{mod}~d)
\end{equation*}%
which is equivalent with 
\begin{equation*}
n\equiv 0(\text{mod}~d)
\end{equation*}%
and 
\begin{equation*}
S_{t}p_{2n+2^{t}-1}-2p_{n}+1\equiv 0(\text{mod}~q).
\end{equation*}%
From the last relation, we have $\widehat{0}=\widehat{2p_{n}-1}$ in $\mathbb{%
Z}_{q}.$ From here, since $n\equiv 0$ (mod $d$), we get $\widehat{2p_{n}-1}=%
\widehat{q-1}$ , which is false. We obtain that there are no idempotent $l-$%
elements in the algebra $A_{t}=\left( \frac{-1,...,-1}{\mathbb{Z}_{q}}%
\right) $.$\Box \smallskip $\smallskip \newline

\textbf{Remark 3.3.} Statement i) from the above Theorem was proved in [Da,
Dr; 70] for the Fibonacci sequence, that means for $l=1$. \newline
\begin{equation*}
\end{equation*}

\bigskip \textbf{4. Examples}%
\begin{equation*}
\end{equation*}%
\qquad

\textbf{Example 4.1. }For $t=2,~$when $\gamma _{1}=-1$ and $\gamma _{2}=q$
is a prime positive integer, it is known, from ([La; 04]), that the
quaternion algebra $\mathbb{H}_{\mathbb{Q}}\left( -1,q\right) $ splits if
and only if $q\equiv 1$ (mod $4$). We wonder how many invertible $l-$%
elements are in the algebra $\mathbb{H}_{\mathbb{Q}}\left( -1,q\right) ?~$
For $q$ \ a prime positive integer, $q\equiv 1$ (mod $4$), we consider the
quaternion algebra $A_{t}=\left( \frac{-1,q}{\mathbb{Q}}\right) =H_{\mathbb{Q%
}}\left( -1,q\right) .$ Let $P_{n}$ be the $n$-th $l-$quaternion. Then, we
have $n\left( P_{n}\right) =p_{2n+1}-q\cdot p_{2n+5},$ therefore the $n$-th $%
l-$quaternion \ $P_{n}$ is invertible for all $n$$\in \mathbb{N}.$ Indeed,
using Proposition 3.1, we have: 
\begin{equation*}
\boldsymbol{n}\left( P_{n}\right) =p_{n}^{2}+p_{n+1}^{2}-q\cdot
p_{n+2}^{2}-q\cdot p_{n+3}^{2}=p_{2n+1}-q\cdot p_{2n+5}.
\end{equation*}%
Remarking that the sequence of $l-$numbers is strictly increasing and $q\geq
1,$ it results that $p_{2n+1}<q\cdot p_{2n+5},$ for all $n$$\in $$\mathbb{N}%
. $ Therefore $\boldsymbol{n}\left( P_{n}\right) \neq 0,$ for all $n$$\in $$%
\mathbb{N}.$ We obtain that $P_{n}$ is invertible, for all $n$$\in $$\mathbb{%
N}.\Box \smallskip $\smallskip

\textbf{Example 4.2}. Even if the octonion algebra $\mathbb{O}\left(
-1,-1,-1\right) =\left( \frac{-1,-1,-1}{\mathbb{R}}\right) $ is division
algebra, for $t=4,$ the real sedenion algebra 
\begin{equation*}
\left( \frac{-1,-1,-1,-1}{\mathbb{R}}\right) ,
\end{equation*}%
with the basis $\{1,e_{1},.....,e_{15}\},$ is not a division algebra. For
example,\newline
$\left( e_{3}+e_{10}\right) \left( e_{6}-e_{15}\right) =0,$(see [Fl; 13] ).
The $l-$element $P_{1}=p_{1}+\overset{15}{\underset{i=1}{\sum }}p_{i+1}\cdot
e_{i}$ has the norm $n\left( P_{1}\right) =M_{2}M_{3}M_{4}p_{17}\neq 0,$
therefore $P_{1}$ is an invertible element.\medskip

Let $q$ be a prime positive integer, $q\geq 3,$ and $\left( \mathbb{Z}%
_{q},+,\cdot \right) $ be the finite field$.$ It is know that the quaternion
algebra $\mathbb{H}_{\mathbb{Z}_{q}}\left( -1,-1\right) $ splits. Therefore,
it has proper zero divisors. In the paper [Sa;17], using many properties of
Fibonacci and Lucas numbers, were determined Fibonacci quaternions and
generalized Fibonacci-Lucas quaternions which are invertible, respectively
zero divisors elements in some quaternion algebras $\mathbb{H}_{\mathbb{Z}%
_{q}}\left( -1,-1\right) .$ Here, for $q=3,$ $q=5,$ respectively $q=7,$ we
get the invertible elements, respectively zero divisors elements and
idempotent elements in the quaternion algebras $\mathbb{H}_{\mathbb{Z}%
_{3}}\left( -1,-1\right) ,$ $\mathbb{H}_{\mathbb{Z}_{5}}\left( -1,-1\right)
, $ $\mathbb{H}_{\mathbb{Z}_{7}}\left( -1,-1\right) .$ A quaternion from $%
\mathbb{H}_{\mathbb{Z}_{q}}\left( -1,-1\right) $ is a zero divisor if and
only if its norm is zero.\medskip\ \smallskip \newline

\textbf{Example 4.3.} Let $q\geq 3$ be a prime positive integer, let $H_{%
\mathbb{Z}_{q}}\left( -1,-1\right) $ be the quaternion algebra and let $%
P_{n} $ be the $n$-th generalized $l-$quaternion. Then, from Proposition
3.1, taking $t=2,~$we have\newline
\begin{equation*}
\boldsymbol{n}\left( P_{n}\right) =\left( l^{2}+2\right) p_{2n+3}.
\end{equation*}

\textbf{Example 4.4.} All Pell quaternions in the quaternion algebra $H_{%
\mathbb{Z}_{3}}\left( -1,-1\right) $ are zero divisors\textit{.} \smallskip
Indeed, we apply Proposition 3.1, that means $l=2,$ and we get $\boldsymbol{n%
}\left( P_{n}\right) =\widehat{0}$ in $\mathbb{Z}_{3},$ for all $n$$\in $$%
\mathbb{N}.\medskip $

\textbf{Example 4.5.} Let $(p_{n})_{n\geq 0}$ be the Pell sequence. Then, we
have: $p_{n}$ $\vdots $ $5$ if and only if $n\equiv 0$ (mod $3$). Indeed, $%
d=3$ is the first number such that $p_{d}$ $\vdots $ $5$ and we apply
Theorem 3.2.\medskip\ \smallskip \newline

\textbf{Example 4.6.} A Pell quaternion is a zero divisor in the quaternion
algebra $H_{\mathbb{Z}_{5}}\left( -1,-1\right) $ if and only if $n\equiv 0$
(mod $3$). Indeed, this statement results from Theorem 3.2.\medskip

\textbf{Example 4.7.} There are no idempotent Pell quaternions in the
quaternion algebra $H_{\mathbb{Z}_{5}}\left( -1,-1\right) $. Indeed, it is a
particular case of Theorem 3.2.\medskip\ \smallskip \newline

\textbf{Example 4.8.} Let $(p_{n})_{n\geq 0}$ be the Pell sequence. Then, we
have: $p_{n}$ $\vdots $ $7$ if and only if $n\equiv 0$ (mod $6$). Obviously, 
$d=6$ is the first number such that $p_{d}$ $\vdots $ $7.$ Therefore, we
apply Theorem 3.2.\medskip

\textbf{Example 4.9.} There are no Pell quaternion zero divisors in the
quaternion algebra $H_{\mathbb{Z}_{7}}\left( -1,-1\right) $. Therefore all
Pell quaternions $P_{n}$ from the quaternion algebra $H_{\mathbb{Z}%
_{7}}\left( -1,-1\right) $ are invertible. It is a consequence of the
Theorem 3.2.%
\begin{equation*}
\end{equation*}

\textbf{Conclusions. }In this paper, we found elements, namely $l-$elements,
with supplementary properties which can help us in the study of algebras $%
A_{t}$ obtained by the Cayley-Dickson process. In Proposition 3.1 and
Theorem 3.2, we found conditions such that these elements are invertible.
These conditions can be very useful in solving equations in algebras $A_{t}$
or to find new identities, identities which can help us in the study of the
properties of these algebras.

\begin{equation*}
\end{equation*}

\textbf{References}\newline
\begin{equation*}
\end{equation*}%
\newline
[Ba, Pr; 09] M. Basu, B. Prasad, \textit{The generalized relations among the
code elements for Fibonacci coding theory}, Chaos, Solitons and Fractals,
41(2009), 2517-2525.\newline
[Ca; 16] P. Catarino, \textit{The Modified Pell and the Modified k-Pell
Quaternions and Octonions}, Adv. Appl. Clifford Algebras, 26(2)(2016),
577-590.\newline
[Ci, Ip; 16] C. B. \c{C}imen. A. Ipek, \textit{On Pell Quaternions and
Pell-Lucas Quaternions}, Adv. Appl. Clifford Algebras, 26(1)(2016), 39-51.%
\newline
[Da, Dr; 70] Daykin, D.E., Dresel, \textit{L.A.G.: Factorization of
Fibonacci numbers}, The Fibonacci Quarterly 8(1), 23\^{a}\euro ``30 (1970). 
\newline
[Fa, Pl; 07(1)] S. Falc\'{o}n, \'{A}. Plaza, \textit{On the Fibonacci} $k$%
\textit{-numbers}, Chaos, Solitons and Fractals, 32(5)(2007), 1615--24.%
\newline
[Fa, Pl; 07(2)] S. Falc\'{o}n, \'{A}. Plaza, \textit{The} $k$\textit{%
-Fibonacci sequence and the Pascal 2-triangle}, Chaos, Solitons and
Fractals, 33(1)(2007), 38--49.\newline
[Fl; 09] C. Flaut, \textit{Divison algebras with dimension }$2^{t}$\textit{, 
}$t\in \mathbb{N}$, Analele \c{S}tiin\c{t}ifice ale Universit\u{a}\c{t}ii
\textquotedblleft Ovidius\textquotedblright\ Constan\c{t}a, Seria
Matematica, 13(2)(2006), 31-38.\newline
[Fl; 13] C. Flaut, \textit{Levels and sublevels of algebras obtained by the
Cayley--Dickson process}, Ann. Mat. Pura Appl., 192(6)(2013), 1099-1114.%
\newline
[Fl, Sa; 15] C. Flaut, D. Savin, \textit{Quaternion Algebras and Generalized
Fibonacci-Lucas Quaternions}, Adv. Appl. Clifford Algebras, 25(4)(2015),
853-862.\newline
[Fl, Sa; 15(1)] C. Flaut, D. Savin, \textit{Some examples of division symbol
algebras of degree 3 and 5}, Carpathian J. Math, 31(2)(2015), 197-204.%
\newline
[Fl, Sa; 17] C. Flaut, D. Savin, \textit{Some remarks regarding} $%
(a,b,x_{0},x_{1})$ \textit{numbers and} $(a,b,x_{0},x_{1})$ \textit{%
quaternions}, submitted, 2017.\newline
[Fl, Sa; 18] C. Flaut, D. Savin, \textit{Some special number sequences
obtained from a difference equation of degree three}, Chaos, Solitons \&
Fractals, 106 (2018), 67-71.\newline
[Fl, St; 09] C. Flaut, M. \c{S}tef\u{a}nescu, \textit{Some equations over
generalized quaternion and octonion division algebras}, Bull. Math. Soc.
Sci. Math. Roumanie, 52(4)(100)(2009), 427-439.\newline
[Ha; 12]\ S. Halici, \textit{On Fibonacci Quaternions}, Adv. in Appl.
Clifford Algebras, 22(2)(2012), 321--327.\newline
[Ha, 15] S.~Halici, On Dual Fibonacci Octonions, Adv. in Appl. Clifford
Algebras, 25(4)(2015), 905--914. \newline
[Ho; 61] A. F. Horadam, \textit{A Generalized Fibonacci Sequence}, Amer.
Math. Monthly, 68 (1961), 455-459.\newline
[Ho; 63]\ A. F. Horadam, \textit{Complex Fibonacci Numbers and Fibonacci
Quaternions}, Amer. Math. Monthly, 70(1963), 289--291.\newline
[La; 04]\ [16] T. Y. Lam, \textit{Introduction to Quadratic Forms over Fields%
}, American Mathematical Society, 2004.\newline
[Mo; 17] G. Cerda-Morales, \textit{On a Generalization for Tribonacci
Quaternions}, Mediterranean Journal of Mathematics, 14(6)2017, 239.\newline
[ Sa; 17] D. Savin, \textit{About Special Elements in Quaternion Algebras
Over Finite Fields}, Advances in Applied Clifford Algebras, June 2017, Vol.
27, Issue 2, 1801-1813.\newline
[ Sa; 18] D. Savin, \textit{Special numbers, special quaternions and special
symbol elements}, accepted for publication in the book Models and Theories
in Social Systems, Springer 2018 (https://arxiv.org/pdf/1712.01941.pdf).%
\newline
[Sc; 66] Schafer, R. D., \textit{An Introduction to Nonassociative Algebras,}
Academic Press, New-York, 1966.\newline
[Si,18] R. da Silva, K.Souza de Oliveira, A.Cunhada, G. Neto, \textit{On a
four-parameter generalization of some special sequences}, Discrete Applied
Mathematics, 243(2018), 154-171\newline
\begin{equation*}
\end{equation*}

\bigskip

Cristina FLAUT

{\small Faculty of Mathematics and Computer Science, Ovidius University,}

{\small Bd. Mamaia 124, 900527, CONSTANTA, ROMANIA}

{\small http://www.univ-ovidius.ro/math/}

{\small e-mail: cflaut@univ-ovidius.ro; cristina\_flaut@yahoo.com}

\medskip \medskip \qquad\ \qquad\ \ 

Diana SAVIN

{\small Faculty of Mathematics and Computer Science, }

{\small Ovidius University, }

{\small Bd. Mamaia 124, 900527, CONSTANTA, ROMANIA }

{\small http://www.univ-ovidius.ro/math/}

{\small e-mail: \ savin.diana@univ-ovidius.ro, \ dianet72@yahoo.com}\bigskip
\bigskip

\end{document}